\documentclass[10pt]{amsart}
\usepackage{amsfonts}
\usepackage{amssymb}
\usepackage{amsthm}
\usepackage{latexsym}
\usepackage{amsmath}
\usepackage[usenames]{color}
\usepackage{soul}

\usepackage[all,dvips]{xy}
\usepackage[dvips]{graphicx}
\newtheorem{teo}{Theorem}
\newtheorem{prop}[teo]{Proposition}

\newtheorem{obs2}[teo]{Remark}
\newtheorem{cor}[teo]{Corollary}
\newtheorem{lema}[teo]{Lemma}

\newenvironment{dem}{\begin{proof}[Proof]}{\end{proof}}

\newenvironment{dem1}{\begin{proof}[Proof of Theorem \ref{hypxn}]}{\end{proof}}
\newenvironment{dem2}{\begin{proof}[Proof of Theorem \ref{biexn}]}{\end{proof}}
\newenvironment{dem3}{\begin{proof}[End of the proof of Theorem \ref{aut}]}{\end{proof}}

\newcommand{\field}[1]{\mathbb{#1}}
\newcommand{\Q}{\field{Q}}
\newcommand{\C}{\field{C}}
\newcommand{\R}{\field{R}}

\newcommand{\Z}{\field{Z}}

\newcommand{\rank}{\rm rank}

\usepackage{amsmath,amssymb,amsfonts,amsthm}

\begin{document}

\title[Bielliptic and Hyperelliptic
modular curves $\boldsymbol{X(N)}$]{Bielliptic and Hyperelliptic
modular curves $\boldsymbol{X(N)}$ and the group
$\boldsymbol{\mathrm{Aut}(X(N))}$.}

%\author{Francesc Bars, Aristides Kontogeorgis and Xavier Xarles}

\author[Francesc Bars]{Francesc Bars}
\address{Francesc Bars\\Departament de Matem\`atiques\\Universitat Aut\`onoma de Barcelona\\08193 Bellaterra, Barcelona, Catalonia}
\email{francesc@mat.uab.cat}

\author[Aristides Kontogeorgis]{Aristides Kontogeorgis}
\address{Aristides Kontogeorgis\\K. Vesri 10 Papagou\\ Athens, GR-15669, Greece}
%\address{Department of Mathematics\\ National and Kapodistrian\\ University of Athens\\ Panepistimioupolis, GR-157 84\\ Athens. Greece}
\email{aristides.kontogeorgis@gmail.com}

\author[Xavier Xarles]{Xavier Xarles}
\address{Xavier Xarles\\Departament de Matem\`atiques\\Universitat Aut\`onoma de Barcelona\\08193 Bellaterra, Barcelona, Catalonia}
\email{xarles@mat.uab.cat}

\thanks{F.Bars and X.Xarles are partially supported by grant MTM2009-10359. A. Kontogeorgis
is  supported by the Thalis Project, Algebraic modeling of
topological and computational structures MIS 380154. }

%\date{23th April, 2012}

\maketitle

\begin{flushright}
    \textsc{Dedicated to the memory of F. Momose}
  \end{flushright}

\begin{abstract} In this paper we determine all modular curves $X(N)$ (with $N\geq 7$) that
are hyperelliptic or bielliptic. We also give a proof that the
automorphism group of $X(N)$ is the group
$\operatorname{PSL}_2(\Z/N\Z)$, therefore it coincides with the
normalizer of $\Gamma(N)$ in $\operatorname{PSL}_2(\mathbb{R})$
modulo $\pm\Gamma(N)$.
\end{abstract}

\begin{section}{Introduction}
In this paper we discuss some basic problems on the modular curves
$X(N)$. By $X(N)$ we mean a geometrically connected curve defined
over $\Q$, which over the complex field $\C$ is given as Riemann
surface by the quotient of $\mathbb{H}^*$ modulo the modular
subgroup
$$\Gamma(N)=\left\{\left(\begin{array}{cc} a&b\\
c&d\\
\end{array}\right)\in\mathrm{SL}_2(\mathbb{Z})\left|\left(\begin{array}{cc} a&b\\
c&d\\
\end{array}\right)\equiv \left(\begin{array}{cc} 1&0\\
0&1\\
\end{array}\right) \pmod N \right.\right\},$$
where as usual we denote by $\mathbb{H}^*$ the union of the upper
half plane $\mathbb{H}$ together with the ``so-called cusps''
$\Q\cup\{\infty\}$. The modular curve $X(N)$ is defined over $\Q$ as
the moduli space parametrizing generalized elliptic curves together
with a full $N$-level structure (see Section 2 for a detailed
discussion).

One of the problems we solve is the determination of the integers
$N\geq 7$ such that the modular curves $X(N)$ are hyperelliptic or
bielliptic. We obtain in Section 4 that none of them are
hyperelliptic (Theorem \ref{hypxn}) and they are bielliptic only for
$N=7$ and $8$ (Theorem \ref{biexn}).

The study of this problem for some other families of modular
curves was initiated by Ogg (see \cite{Og} and \cite{ba}) with the
case of the curves $X_0(N)$, and followed by the modular curves
$X_1(N)$ in \cite{Mes} and \cite{jeki} and the curves $X_1(N,M)$
in \cite{IsMo} and \cite{jeki2}.

In the last section we apply this result to study the finiteness of
quadratic points of $X(N)$. In particular we can prove that the set
of quadratic points of $X(N)$ over the cyclotomic field
$\Q(\zeta_N)$ is always  finite for $N>6$, $\zeta_N$ denotes, as
usual, a primitive $N$-root of unity.

In Section 3 we consider another important issue concerning the
curve $X(N)$: the explicit determination of its automorphism group
over $\C$, that we denote by $\mathrm{Aut}(X(N))$.

Recall that, for a modular curve $X$ of genus greater than one and
with modular group $\Gamma\leq \operatorname{SL}_2(\Z)$, the
quotient of the normalizer of $\Gamma$ in
$\operatorname{PSL}_2(\mathbb{R})$ by $\pm\Gamma$ always gives a
subgroup of $\operatorname{Aut\,}(X)$. We denote this subgroup by
$\mathrm{Norm}(\Gamma)/\pm\Gamma$. It is a quite difficult problem
to determine when $\mathrm{Norm}(\Gamma)/\pm\Gamma$ coincides with
the full group of automorphisms of the corresponding modular curve
$X$, of genus greater than one. An automorphism $v\in
\mathrm{Aut}(X)\setminus (\mathrm{Norm}(\Gamma)/\pm\Gamma)$ is
called exceptional.

Kenku and Momose in \cite{KeMo} determined the full automorphism
group for $X_0(N)$ with $N\neq 63$ and Elkies in \cite{Elkies}
obtained $\mathrm{Aut}(X_0(63))$, and finally Harrison in \cite{Har}
corrected the Kenku-Momose statement for
$\mathrm{Aut}(X_0(108))$\footnote{We would like to mention that this
correction does not affect the results in \cite{ba}.}. In
particular, there are exceptional automorphisms for $X_0(N)$ only
for $N= 37,63$ and $108$. For the modular curve $X_1(N)$ and $N$
square-free, Momose proved in \cite{Mo} that there are no
exceptional automorphism.

Finally, let us explain briefly the history concerning
$\mathrm{Aut}(X(N))$. J.P. Serre in a letter to B. Mazur \cite{se}
computed that the automorphism groups of the modular curves $X(p)$
for $p$ prime $p\geq 7$, are isomorphic to the simple groups
$\mathrm{PSL}_2(\Z/p\Z)$. Back in 1997, in a conference held in Sant
Feliu de Gu\'ixols, the second author met G. Cornelissen, who wanted
to compute the automorphism group of the Drinfeld modular curves
\cite[sec. 10]{CKK} and asked if there is a generalization for
composite $N$. After finishing the computation, the second author
communicated the generalization to J.P. Serre, who answered that the
theorem he proved in the letter to Mazur and the generalization for
composite $N$ should be known to the pioneers of modular forms; but
we were not able to find a reference in the literature. Since there
is new interest (see \cite{Gur},\cite{Ritz}) on the automorphisms of
the modular curves $X(N)$ and we are not aware of any reference for
this computation, we believe that writing down a proof might be
useful to the mathematical community.

\end{section}

\begin{section}{Preliminary results on the curve $X(N) $}

The (non-complete) Riemann surface $Y(N)(\C)$ is constructed as
the quotient of the upper-half plane $\mathbb{H}$ modulo the
modular subgroup $\Gamma(N)$. The set $Y(N)(\C)$ parametrizes the
pairs $(E,(P_1,P_2))$, where $E$ is an elliptic curve defined over
$\C$ and $P_1$ and $P_2$ are points of order exactly $N$ in
$E(\C)$, which generate the subgroup of $N$-torsion points and
verify that $e(P_1,P_2)=\exp(\frac {2\pi i}{N})$, where $e$
denotes the Weil pairing.

This interpretation can be used to give a model of the modular curve
$Y(N)$ (and its completion $X(N)$) over other fields of
characteristic not dividing $N$ (or general schemes over
$\mbox{Spec}(\mathbb Z[1/N]$). We have two options: either we ignore
the last condition on the Weil paring, obtaining a non geometrically
connected curve, or we modify the moduli problem introducing in some
way the Weil pairing. The first option is essentially the one taken
by Deligne and Rappoport in \cite{DR}, and also by Katz and Mazur in
\cite{KM}. We consider here the second option, following for example
Elkies in \cite{Elkies2}, Section 4. Over fields containing all the
$N$-roots of unity $\mu_N$ the second curve is isomorphic to a
connected component of the first.

Thus, we take the full modular curve $Y(N)$ (for $N>2$) as the
(geometrically connected) curve which over any field $K$ (of
characteristic not dividing $N$) parametrizes pairs $(E,\phi)$,
where $E$ is an elliptic curve over $K$ and $\phi$ is a
Weil-equivariant isomorphism of group schemes between $E[N]$, the
kernel of the multiplication by $N$ on $E$, and $\Z/N\Z\times
\mu_N$. This means that $$\langle \ ,\ \rangle\circ (\phi\times
\phi)=e$$ where $e:E[N]\times E[N] \to \mu_N$ is the Weil Paring and
$\langle\ , \ \rangle$ is the natural (symplectic) self-pairing of
$\Z/N\Z\times \mu_N$ given by
$$ \langle (m,\xi),(n,\eta) \rangle:=\xi^n\eta^{-m}.$$
The map $\phi$ is called the $N$-level structure. We denote by
$X(N)$ the completion of $Y(N)$; it also has a moduli interpretation
like $Y(N)$ by allowing
 generalized elliptic curves.  For the
cases $N=1$ and $2$ one takes the usual coarse moduli space (in both
cases isomorphic to the projective line).

There are other options one can take to get a model of the curve
$X(N)$; for example, one can take a fixed elliptic curve
$\tilde{E}$, and consider the $N$-level structures $\phi$ given as
Weil-equivariant isomorphisms between $E[N]$ and $\tilde{E}[N]$.
One gets a twisted form of $X(N)$, usually denoted by
$X_{\tilde{E}}(N)$ (see for example \cite{HK}).

Recall that the curve $X(1)$ is isomorphic (via the $j$-function) to
the projective line $\mathbb{P}^1$. The canonical cover
$X(N)\rightarrow X(1)$ that forgets the $N$-level structure is
Galois over any field containing all the $N$-roots of unity, and
with Galois group $\mathrm{PSL}_2(\mathbb{Z}/N \mathbb{Z})$. The
degree of this cover is equal to
$$
\delta _N:=\left\{
\begin{array}{ll}
N^3/2\prod_{p\mid N}(1-p^{-2}) & \text{if }N>2 \\
6 & \text{if }N=2
\end{array}
\right..
$$
Moreover the genus $g_N$ of  $X(N)$ is equal  to \cite[p. 23]{Shi}
\begin{equation}
g_N=1+\delta _N\frac{N-6}{12N}.  \label{genus}
\end{equation}
We see that the curves  $X(2),X(3),X(4)$ and $X(5)$ are rational,
while the curve  $X(6)$ is elliptic. For all the other values $N\geq
7$ the curves $X(N)$ have genus $>1$.

We now want to relate the curve $X(N)$ to some other modular
curves. First, observe that we have natural ``forgetful" maps
$f_1:X(N) \to X_1(N) $ given, in the moduli interpretation, by
sending the pair $(E,\phi)$ to the pair $(E,\phi^{-1}((1,1)))$,
since $\phi^{-1}((1,1))$ is a point of exact order $N$. Thus, we
have maps $f_0:X(N) \to X_0(N) $ obtained composing the map $f_1$
with the forgetful map $\varrho:X_1(N)\to X_0(N)$. There is also
another independent map $f_0':X(N) \to X_0(N)$, which  can be
defined  in terms of the moduli interpretation, as the map sending
the pair $(E,\phi)$ to the pair $(E,\phi^{-1}(\{0\}\times
\mu_N))$. If we see the curve $X(N)$ as the compactified quotient
of $\mathbb{H}$ by a discrete subgroup, then we can interpret
these maps  $f_1$ and $f_0'$ as the quotient maps of $X(N)$ by the
subgroups
$$\Gamma_1(N)=\left\{ \left(\begin{array}{cc}
1&*\\
0&1\\
\end{array} \right)\right\}\mbox{ and  } \Gamma^0(N)=\left\{ \left(\begin{array}{cc}
*&0\\
*&*\\
\end{array} \right)\right\} \subseteq \mathrm{PSL}_2(\Z/N\Z),$$ respectively.

Over a field containing a primitive $N$-root of unity $\zeta_N$,
there is a map $f_1':X(N)\to X_1(N)$, which depends on $\zeta_N$,
given by assigning to the pair $(E,\phi)$, in the notation above,
the $N$-torsion point $\phi^{-1}((0,\zeta_N))$. The map $f_0'$ can
be factored as $f_0'=\varrho\circ f_1'$.

We now recall a construction of natural maps from  $X_1(N^2)$ to
$X(N)$ and from $X(N)$ to $X_0(N^2)$, for which we do not know a
precise reference (see, however, Section 11.3.5. in \cite{KM}, for
the second morphism in the case $N=p^n$, $p$ a prime).

\begin{lema}\label{X0N2} Let $N\geq3$ be an integer. Then there exist
morphisms of curves $\pi_1:X_1(N^2)\to X(N)$ of degree $N$ and
$\pi_0:X(N)\to X_0(N^2)$ of degree $\varphi(N)/2$ defined over
$\Q$, such that the composition $\pi_0\circ \pi_1:X_1(N^2)\to
X_0(N^2)$ is the natural forgetful map. Moreover, the maps make
the following diagram commutative:
$$
\xymatrix@R=0.8pc@C=0.9pc{
&        &        &      \ar@/_8mm/[ddddlll] X_1(N^2) \ar@{->}[d]_{\pi_1}   &         &       &     \\
&        &        &      \ar@/_4mm/[dddlll]_{f_1} X(N) \ar@{->}[ddl]_{\pi_0}    \ar@/^10mm/[ddddrrr]^{f'_0} &         &       &     \\
&        &        &                                                                 &         &       &     \\
&        &         \ar@/_1mm/[dd]     X_0(N^2)      \ar@{<->}[rr]^{\omega_{N^2}}  &  &    X_0(N^2)      \ar@/^1mm/[dd]    &    &   & \\
X_1(N) \ar@{->}[d] &       &            &       &           &                          &   \\
 X_0(N)\ar@{<->}[rr]^{\omega_N} \ar@{->}[ddrrr]  &    & X_0(N)             &   & X_0(N) \ar@{<->}[rr]^{\omega_N}   &           & \ar@{->}[ddlll] X_0(N)
 \\ \\
 &    &            &  X(1)     &           & & &
\\}
$$
where $\omega_N$ and $\omega_{N^2}$ denote the Atkin-Lehner
involutions, and the maps without name are the usual projection maps
given by the forgetful maps.
\end{lema}

\begin{dem}

We will construct the maps from $X_1(N^2)$ to $X_0(N)$ and from
$X(N)$ to $X_0(N^2)$ in two equivalent ways. First, over the complex
numbers, the map is deduced by observing that
$$\Gamma_1(N^2) \leq U^{-1} \Gamma(N) U \leq \Gamma_0(N^2),$$
where
$$U=\left(\begin{array}{cc} 1&0\\ 0&1/N
\end{array}\right).$$
These maps can be defined over $\Q$ (or any field with
characteristic prime to $N$) by using the moduli interpretation.
First, the map from $X(N)$ to $X_0(N^2)$ can be described on $Y(N)$
by sending the point of $Y(N)$ given by an elliptic curve $E$ and
the $N$-level structure $\phi:E[N]\to \Z/N\Z\times\mu_N$ to the
$N^2$-cyclic isogeny obtained composing the dual of the $N$-isogeny
$E\to E/F_1$ with the $N$-isogeny $E\to E/F_2$, where we consider
the subgroups $F_1:=\phi^{-1}(\Z/N\Z\times \{1\})$ and
$F_2:=\phi^{-1}(\{0\}\times \mu_N$).

The morphism $\pi_0$ can also be interpreted as the natural map from
$X(N)$ to $X(N)/C$, where $C$ is the full Cartan subgroup of
$\mathrm{PSL}_2(\Z/N\Z)$ (formed by the diagonal matrices).

The map from $X_1(N^2)$ to $X(N)$ can be analogously described in
the moduli interpretation for the points in $Y_1(N^2)$ over a
field $K$, given as pairs $(E,P)$ where $E$ is an elliptic curve
over $K$ and $P$ is a point of exact order $N^2$: consider the
point $Q:=[N]P$, which has order $N$, and the elliptic curve
$E':=E/ \langle Q \rangle$. Then $E'$ has two natural cyclic
isogenies of degree $N$, the quotient $E'\to E/ \langle P \rangle$
and the dual isogeny of $E\to E'$. The kernel $F_1$ of the first
map is canonically  isomorphic to $\Z/N\Z$, i.e. $F_1\cong
\Z/N\Z$, where the isomorphism is given by the point $P$. We
denote by $F_2$ the kernel of the second isogeny $E'\to E$, dual
of $E\to E'$. Then the two subgroups $F_1$ and $F_2$ have zero
intersection and hence there must be a canonical isomorphism
$F_2\cong \mu_N$ given by the Weil pairing. Therefore we have an
Weil-equivariant isomorphism $\phi:E'[N]=F_1\oplus F_2\cong
\Z/N\Z\times \mu_N$.

The commutativity of the diagram is clear from the definition of
the maps via the moduli interpretation of the curves. Recall that
the natural projection map from $X_0(N^2)$ to $X_0(N)$ sends a
non-cuspidal point $(E,\varphi)$ of $Y_0(N^2)$ to the point
$(E,\varphi_1)$, where $\varphi=\varphi_2\circ \varphi_1$ is the
decomposition of the degree cyclic $N^2$ isogeny $\varphi:E\to E'$
as composition of two cyclic degree $N$ isogenies, and that the
Atkin-Lehner involution sends an isogeny to is dual.

Finally, the assertions on the degrees are easy over $\C$, taking
into account that the subgroup $\Gamma_0(N^2)$ contains
$-\mathrm{Id}_2$, but $\Gamma(N)$ (and $\Gamma_1(N^2)$) do not.
\end{dem}

\begin{obs2} This lemma implies that, for $N=3$, $4$ and $6$, the moduli curves
$X(N)$ and $X_0(N^2)$ are identical over $\Q(\zeta_N)$. This is
analogous to the case of the curves $X_1(N)$ and $X_0(N)$ for
$N=3$, $4$ and $6$. Note that this does not imply that given an
elliptic curve $E$ over a field $K$ and a cyclic subgroup scheme
$F$ of order $3$ defined over $K$, then $F$ contains a point of
order $3$; but that there exists a (unique) quadratic twist $E'$
of $E$ such that the corresponding subgroup scheme $F'$ of $E'$
contains a point of order $3$. Equivalently, there exists a point
$P$ in $F$, defined over a (quadratic) extension $L$, such that
the pair $\{P,-P\}$ is defined over $K$. The same is true for
$N=4$ and $6$, and, in general, for the elliptic curves whose
$j$-invariant is in the image of the map $Y_1(N)(K)\to
Y(1)(K)\overset{\scriptstyle j}{\to} K$.
\end{obs2}

\begin{cor} The curve $X(N)$ is isomorphic over $\Q$ to the fiber
product of $X_1(N)$ and $X_0(N^2)$ over $X_0(N)$, with respect to
the natural map $X_1(N)\to X_0(N)$ and the map $X_0(N^2)\to
X_0(N)\overset{\scriptstyle \omega_N}{\to}X_0(N)$ given by the
composition of the natural map with the Atkin-Lehner involution
$\omega_N$.
\end{cor}

\begin{dem} From the previous lemma and the universal property of
the fiber product we have a natural map from $X(N)$ to the fiber
product. In order to show it is an isomorphism we will prove they
both parametrize the same moduli problem. The moduli problem
parametrized by the fiber product is easily seen to be the
triplets $(E,P,\varphi)$ where $E$ is an elliptic curve, $P$ is a
point of order exactly $N$, $\varphi:E''\to E/\langle P\rangle$ is
a degree $N^2$ cyclic isogeny such that
$\varphi=\varphi_2\circ\varphi_1$, where $\varphi_1:E\to E/\langle
P\rangle$ is the quotient isogeny. Now, the kernel of the dual of
$\varphi_2$ is a subgroup scheme $F$ of order $N$ in $E$. From the
condition $E/F\to E\to E/\langle P\rangle$ being a cyclic isogeny
of degree $N^2$, it is deduced that the subgroups $F$ and $\langle
P\rangle$ have zero intersection. Hence $E[N]\cong F\times \langle
P\rangle\cong F\times \Z/N\Z$. The Weil pairing implies then that
$F\cong\mu_N$ and that this isomorphism is compatible with the
Weil pairing.
\end{dem}

\end{section}

\begin{section}{The automorphism group of $X(N)$}

Recall that the curves $X(N)$ have genus greater than two if $N\ge
7$, and their automorphism groups are bounded by Hurwitz bound:
\begin{equation}
|\mathrm{Aut}(X(N)|\leq 84(g_N-1).  \label{Hur}
\end{equation}
It is also known that exactly three points of $X(1)$ are ramified
in the cover  $X(N)\rightarrow X(1)$, namely $j(i),j(\omega
),j(\infty )$, with ramification indices $2,3$ and   $N$,
respectively ($j$ denotes the natural $j$-invariant isomorphism
between $\mathbb{P}^1$ and $X(1)$). The main result of this
section is the following:
\begin{teo}\label{aut}
 The automorphism group of $X(N)$ over $\C$ for values $N$ such that $g_N\geq 2$ equals
 $\mathrm{PSL}_2(\mathbb{Z}/N\mathbb{Z})$.
\end{teo}

We will proof the theorem in several steps.

\begin{lema} \label{normal}
If   $\mathrm{\mathrm{PSL}}(2,\mathbb{Z}/N\mathbb{Z})\lhd
\mathrm{Aut}X(N)$ then
$\mathrm{\mathrm{PSL}}(2,\mathbb{Z}/N\mathbb{Z})= \mathrm{Aut}X(N)$.
\end{lema}
\begin{dem}
Since $\mathrm{\mathrm{PSL}}(2,\mathbb{Z}/N\mathbb{Z})\lhd
\mathrm{Aut}X(N)$, we can restrict  automorphisms in
$\mathrm{Aut}X(N)$
 to  automorphisms of $X(1)\cong \mathbb{P}^1$ and these automorphisms should
fix the three ramification points. Therefore the restriction  is
the identity.
\end{dem}

Let  $m$ be the index of
$\mathrm{PSL}_2(\mathbb{Z}/N\mathbb{Z})=\mathrm{Gal}(X(N)/X(1))$ in
  $\mathrm{Aut}X(N)$.
The equation for the genus (\ref{genus}) for  $N\neq 2$ can be
written as
\begin{equation}
84(g_N-1)=|\mathrm{PSL}_2(\mathbb{Z}/N\mathbb{Z})|\left(
7-\frac{42}N\right)  \label{genus2}
\end{equation}
and this combined with (\ref{Hur}) gives the following bounds for
the index $m$:
\begin{equation}
\begin{array}{ccc}
m\leq 2 & \mbox{for} & 7\leq N<11 \\
m\leq 3 & \mbox{for} & 11\leq N<14 \\
m\leq 4 & \mbox{for} & 14\leq N<21 \\
m<7 & \mbox{for} & 21\leq N
\end{array}
.  \label{index-bounds}
\end{equation}
Therefore, for  $7\leq N<11$ we have $\mathrm{Aut}X(N)\cong
\mathrm{PSL}_2(\mathbb{Z}/N\mathbb{Z})$ by lemma \ref{normal}.

The following lema is elementary.
\begin{lema} \label{beta}
Consider the coset decomposition
\[
 \mathrm{Aut}(X(N))=
 a_1\mathrm{PSL}_2(\mathbb{Z}/N\mathbb{Z})\cup\cdots \cup a_m\mathrm{PSL}_2(\mathbb{Z}/N\mathbb{Z})
\]
and define the  representation
\[
\beta :\mathrm{PSL}_2(\mathbb{Z}/N\mathbb{Z})\longrightarrow S_m
\]
by sending
\[\sigma \mapsto \{\sigma a_1\mathrm{PSL}_2({\mathbb{Z}}/{N\mathbb{Z}}),\sigma
a_2\mathrm{PSL}_2(\mathbb{Z}/N\mathbb{Z}),...,\sigma
a_m\mathrm{PSL}_2(\mathbb{Z}/N\mathbb{Z})\}.\] Then
$\mathrm{PSL}_2(\mathbb{Z}/N\mathbb{Z})\lhd \mathrm{Aut}X(N)$, if
and only if  $\beta $ is the trivial homomorphism.
\end{lema}

\begin{lema}
 If $N=p$ is prime, $p\geq 7$ then $\beta=1$.
\end{lema}
\begin{dem}
Since  $\mathrm{PSL}_2(p)$ is simple we have  $\mathrm{ker}\beta$ is
either  $\mathrm{PSL}_2(p)$ or   $\{1\}$. The last case is
impossible since there are no elements of order  $p$ in   $S_m$, for
 $m\leq 6$.
\end{dem}
Let us now consider the curves  $X(p^e)$, where   $p$ is prime,
$p\geq 7$.
\begin{lema}
For $X(p^e)$ with $p\geq 7$ we have $\mathrm{Aut}
X(p^e)=\mathrm{PSL}_2(p^e)$.
\end{lema}
\begin{dem}
We will prove that $\beta=1$ for the map $\beta$ defined in lemma
\ref{beta}. We consider the following tower of covers
\[
 \xymatrix{X(p^e)
 \ar[d]^H \ar@/_2pc/_{\mathrm{PSL}_2(\mathbb{Z}/p^e\mathbb{Z})}[dd]\\
 X(p) \ar[d]^{\mathrm{PSL}_2(p)} \\
 X(1)
 }
\]
Consider $H:=\mathrm{Gal}(X(p^e)/X(p))$; then   $|H|=p^{3(e-1)}$,
and, since $p\geq 7$, we have $H<\ker \beta $. Therefore we can
define the homomorphism  $\tilde{\beta}$  so that the following
diagram is commutative
\[
 \xymatrix{ \mathrm{PSL}_2(\mathbb{Z}/p^e\mathbb{Z})\ar[d]  \ar^{\beta}[rr]& & S_m \\
  \mathrm{PSL}_2(p) \ar_{\tilde{\beta}}[rru] & &
}
 \]
Again, since $\mathrm{PSL}_2(p)$ is simple, we obtain
$\tilde{\beta}=1$ and the same holds for  $\beta$.
 \end{dem}
\begin{cor}
Let  $N$ be a composite integer prime to   $2,3,5$. Then
$\mathrm{Aut}X(N)=\mathrm{PSL}_2(\mathbb{Z}/N\mathbb{Z})$.
\end{cor}
\begin{dem}
The homomorphism   $\beta $ is trivial in this case as well, since
\[
\mathrm{PSL}_2(\mathbb{Z}/N\mathbb{Z})\cong \bigoplus_{i=1}^s \mathrm{PSL}_2(\mathbb{Z}/p_i^{a_i}\mathbb{Z)}%
,
\]
where $N=\prod_{i=1}^sp_i^{a_i}$ is the decomposition of  $N$ in
primes.
\end{dem}

\begin{dem3}
In order to study the case for general $N$ we will need better
bounds for the index
\[m:=[\mathrm{Aut}X(N):\mathrm{PSL}_2(\mathbb{Z}/N\mathbb{Z}].\]
We consider the tower of covers
\[
 \xymatrix{
 X(N) \ar[d]^{\mathrm{PSL}_2(\mathbb{Z}/N\mathbb{Z})} \ar@/_2pc/_{\mathrm{Aut}(X(N))}[dd]&
 \ar@{-}[d]_2 & \ar@{-}_3[d]& \ar@{-}_N[d] \\
 X(1) \ar[d] & j(i) & j(\omega) & j(\infty) \\
 X(N)^{\mathrm{Aut}X(N)}
 }
\]
Observe that if  $\mathrm{PSL}_2(\mathbb{Z}/N\mathbb{Z})$ is not a
normal subgroup of $\mathrm{Aut}X(N) $ then the cover  $X(1)\cong
\mathbb{P}^1\rightarrow X(N)^{\mathrm{Aut}X(N)}$ is not Galois. From
the proof of the Hurwitz bound (\ref{Hur}) for the size of the
automorphism group of an algebraic curve, found in \cite[p.
260]{F-K}, we see that if the number  $r$ of points of $X(1)$
ramified in the cover $X(N)\rightarrow X(N)^{\mathrm{Aut}X(N)}$ is
$r>3$, then Hurwitz's bound is improved to
\[
|\mathrm{Aut}X(N)|\leq 12(g_N-1).
\]
This proves that  $m\leq 1$, so
$\mathrm{PSL}_2(\mathbb{Z}/N\mathbb{Z)}\lhd \mathrm{Aut}X(N),$ a
contradiction.

Therefore the number of ramified points is  $r=3$. Now Hurwitz's
bound for $X(N)\longrightarrow X(N)^{\mathrm{Aut}X(N)}$ gives
\begin{equation}
2(g_N-1)=|\mathrm{Aut}X(N)|\left( 1-\frac 1{\nu _1}+1-\frac 1{\nu
_2}+1-\frac 1{\nu _3}-2\right) ,  \label{Riem-Hur}
\end{equation}
where  $\nu _i$ are the ramification indices of the ramified points
of the cover $X(N)\rightarrow X(N)^{\mathrm{Aut}X(N)}$. We
distinguish the following cases:

{\bf Case 1}  The three points  $j(i),j(\omega ),j(\infty )$
restrict to different points $p_1,p_2,p_3$ with ramification indices
$e(j(i)/p_1)=\kappa ,e(j(\omega )/p_2)=\lambda ,e(j(\infty
)/p_3)=\mu$.  Equation  (\ref{Riem-Hur}) in this case is written
\begin{eqnarray*}
2(g_N-1) &= &|\mathrm{Aut}X(N)|\left( 1-\frac 1{2\kappa }+1-\frac
1{3\lambda }+1-\frac
1{N\mu }-2\right) \geq \\
&\geq & |\mathrm{Aut}X(N)|\left( 1-1/2+1-1/3+1-1/N-2\right) \geq  \\
&\geq & |\mathrm{Aut}X(N)|\left( 1/6-1/N\right)
\end{eqnarray*}
which in turn gives the desired result
\[
|\mathrm{Aut}(X)|\leq \frac{12N}{N-6}(g_N-1)=\delta _N.
\]

{\bf Case 2.}  Some of the three points $j(i),j(\omega ),j(\infty
)$ restrict to the same point $X(N)^{\mathrm{Aut}X(N)}.$  We will
consider the case  $N\geq 11$. First, let us see that the points
$j(i)$ and $j(\infty )$ could not restrict to the same point of
 $X(N)^{\mathrm{Aut}X(N)}.$ Since the cover  $X(N)\rightarrow
X(N)^{\mathrm{Aut}X(N)}$ is  Galois we should have $2\kappa =N\mu$
(using the notations as in the Case 1). But the degree of the
cover $X(1)\rightarrow X(N)^{\mathrm{Aut}X(N)}$ is at most $m\leq
6,$ so $ \kappa \leq m \leq 6$ and $\mu =1,$ and this means that
$j(i)$ and $j(\infty )$ could not restrict to the same point,
unless $N\leq 12.$ But if $N\leq 12$ then $\kappa\leq m\leq 3,$ so
$N\leq 6,$ which contradicts $N\geq 11.$ Using the same argument
we can show that the points $j(\omega )$ and  $j(\infty )$
restrict to different points of $X(N)^{\mathrm{Aut}X(N)}.$

Hence, we can suppose that only the points $j(i)$ and $j(\omega )$
restrict to the same point of $X(N)^{\mathrm{Aut}X(N)}$. Thus,
there should be another point $p$ of $X(N)^{\mathrm{Aut}X(N)}$
which ramifies only in the cover $X(1)\rightarrow
X(N)^{\mathrm{Aut}X(N)}$ with ramification index $2\leq \nu \leq
6.$ The Hurwitz bound implies
\[
2(g_N-1)=|\mathrm{Aut}X(N)|\left( 1-\frac 1{6\psi }+1-\frac 1{\phi
N}+1-\frac 1\nu -2\right) \stackrel{\nu \geq 2,\psi =1\text{ or
}2}{\geq }
\]
\[
|\mathrm{Aut}X(N)|\left( \frac 13-\frac 1N\right) \stackrel{N\geq 11}{\geq }%
|\mathrm{Aut}X(N)|\left( \frac 13-\frac 1{11}\right)
\]
which gives
\[
|\mathrm{Aut}X(N)|\leq 33/4\left( g_N-1\right)
\]
and in turn gives the desired result  $m\leq 1$. \end{dem3}

%\begin{obs2}
%It is interesting to point out that since  $X(2),X(4),X(3),X(5),$
%are rational curves, the classification of finite subgroups of the
%rational function field implies the well known result
%\[
%\mathrm{PSL}_2(2)\cong
%{A}_3,\mathrm{PSL}_2(\mathbb{Z}/4\mathbb{Z)}\cong
%S_4,\mathrm{PSL}_2(3)\cong A_4,\mathrm{PSL}_2(5)\cong A_5
%\]
%\end{obs2}

Recall that $\operatorname{Aut\,}(\mathbb{H})$ is isomorphic to
$\operatorname{PSL}_2(\R)$, and $\Gamma(N)$ is torsion-free if
$N\geq 5$, thus the automorphism group of
$Y(N)=\mathbb{H}/\Gamma(N)$ is the quotient of the normalizer of
$\Gamma(N)$ in $\operatorname{PSL}_2(\R)$ by $\pm\Gamma(N)$,

\begin{cor} For $N\geq 7$  we have $\operatorname{Aut\,}(Y(N))\cong
\mathrm{SL}_2(\Z/N\Z)/\pm1$ and the order of the group of automorphisms of
$Y(N)$ is given by
$$\frac{1}{2}N\varphi(N)\psi(N)$$
where $\varphi(N):=N\prod_{p|N}(1-p^{-1})$ and
$\psi(N):=N\prod_{p|N}(1+p^{-1})$ with $p$ prime.

In particular the normalizer of $\Gamma(N)$ in
$\operatorname{PSL}_2(\R)$ is given by $\operatorname{PSL}_2(\Z)$
and $\mathrm{Norm}(\Gamma(N))/\pm\Gamma(N)\cong\mathrm{PSL}_2(\Z/N\Z)$.
\end{cor}
\begin{dem}
Clearly for $N\geq 5$, $\Gamma(N)\leq \operatorname{PSL}_2(\Z)\leq
\mathrm{Norm}(\Gamma(N))\leq \operatorname{PSL}_2(\R)$. Thus
$$\operatorname{PSL}_2(\Z)/\pm\Gamma(N)\leq
\mathrm{Norm}(\Gamma(N))/\pm\Gamma(N)=\mathrm{Aut}(Y(N))\leq\mathrm{Aut}(X(N))$$ but
$\operatorname{PSL}_2(\Z)/\pm\Gamma(N)$ is isomorphic to $\mathrm{Aut}(X(N))$
for $N\geq 7$ therefore the result.
\end{dem}
\begin{obs2} Following the proof for computing $\mathrm{Norm}(\Gamma_1(N))$ in \cite{kiko}
one can easily deduce that, for $N\geq 5$, $\mathrm{Aut}(Y(N))\cong
\mathrm{PSL}_2(\Z/N\Z)$ and the normalizer of $\Gamma(N)$ in
$\operatorname{PSL}_2(\R)$ is $\operatorname{PSL}_2(\Z)$.
\end{obs2}

\end{section}

\begin{section}{Hyperelliptic and bielliptic modular curves $X(N)$}

Recall that a non-singular projective curve $C$ of genus $g_C>1$
over an algebraically closed field of characteristic zero is
hyperelliptic if it has an involution $v\in
\operatorname{Aut\,}(C)$, called hyperelliptic involution, which
fixes $2g_C+2$ points (see, for example, \cite[\S1]{Sch}). This
involution $v$ is unique if $g_C\geq 2$. Similarly, the curve $C$
is bielliptic if it has an involution $w\in
\operatorname{Aut\,}(C)$, named bielliptic, which fix $2g_C-2$
points. This involution is unique if $g_C\geq 6$.

In this section we want to determine for exactly which integers
$N\ge 7$ the curve $X(N)$ is hyperelliptic or bielliptic over $\C$.
Since $X(N)$ is naturally isomorphic over the cyclotomic field
$\Q(\zeta_N)$ to the curve $X_1(N,N)$, these results can also be
deduced from the results by Ishii-Momose in \cite{IsMo} in the
hyperelliptic case, and by Jeon-Kim in \cite{jeki2} in the
bielliptic case \footnote{Some results in \cite{IsMo} use the no
existence of exceptional automorphisms for intermediate modular
curves \cite{Mo}, but Andreas Schweizer communicated to us that
there are exceptional automorphisms in some intermediate curve (see
the forthcoming work \cite{jeki3}). This correction does not affect
the result on $X_1(N,N)$ in \cite{IsMo} and \cite{jeki2}, but here
we present a proof without using any of the results stated in
\cite{Mo}.}. Here we present a distinct and direct proof.

\begin{teo}\label{hypxn} For $N\geq 7$ the modular curve $X(N)$ is not hyperelliptic.
\end{teo}

\begin{teo}\label{biexn} For $N\geq7$ the modular curve $X(N)$ is bielliptic only when $N=7$
or $N=8$.
\end{teo}

Before we proceed to the proof of the theorems, we collect some
results we will use. Observe first that, given a morphism of
non-singular projective curves
$$\phi:X\rightarrow Y$$  which is a Galois cover (in the sense that it is given by
a quotient map of the form $X\to X/H$, for $H$ a subgroup of the
group of automorphisms of $X$), and given $\nu$ an involution on
$X$, if $\nu$ satisfies $\nu H=H\nu$, then either $\nu$ induces, by
$\phi$, an element of the Galois group $H$ of the cover, or it
induces an involution on $Y$.

\begin{lema}\label{lema3} Consider a Galois cover $\phi:X\rightarrow Y$ of degree $d$
between two non-singular projective curves of genus $g_X\geq 2$ and
$g_Y$ respectively. Suppose that $g_Y\geq 2$ or $d$ is odd.
\begin{enumerate}
\item  Suppose that $2g_X+2>d(2g_Y+2)$. Then, $X$
is not hyperelliptic.

\item Denote by $n_{\iota}$ the number of fixed points of an
involution $\iota$ of $Y$. Suppose $2g_X-2>dn_{\iota}$ for any
involution $\iota$ on $Y$. Then, if $g_X\geq 6$, $X$ is not
bielliptic.

\item Suppose $2g_X-2>d(2g_Y+2)$. Then, if $g_X\geq 6$, $X$ is not
bielliptic.
\end{enumerate}
\end{lema}

\begin{dem} If $v$, a hyperelliptic or bielliptic involution, is
in the group of the Galois cover $\phi$, then we have the
following factorization of $\phi$
$$X\rightarrow X/\langle v \rangle \rightarrow Y,$$ which is
impossible if $d$ is odd, since $X\rightarrow X/\langle v \rangle$
has degree 2, and also if $g_Y\geq 2$, since $X/\langle v \rangle$
has genus $\leq 1$.

Suppose now that $X$ has a hyperelliptic or bielliptic involution
$v$, which induces an involution $\tilde{v}$ on $Y$. Then the
involution $v$ can have fixed points only among the points lying
above the fixed points by $\tilde{v}$ of $Y$, and hence the map $v$
has at most $d n_{\tilde{v}}$ fixed points, where $n_{\tilde{v}}$
denotes the number of fixed points of $\tilde{v}$ on $Y$. By the
Hurwitz's formula, the involution $v$ must have $2g_X+2$ fixed
points in the hyperelliptic case, or $2g_X-2$ fixed points in the
bielliptic case. We get the result under our hypothesis, since
hyperelliptic involutions and bielliptic involutions on $X$ are
(unique and) in the center of $\mathrm{Aut}(X)$ (see
\cite[Proposition 1.2]{Sch}).
\end{dem}

The following lemmas can be easily proved over $\C$ by observing
that both curves attain the maximal order of the group of
automorphisms for their genus. Recall that from the main result in
Section 3 we have $\mathrm{SL}_2(\Z/N\Z)/\pm1\cong
\operatorname{Aut\,}(X(N))$. Now, the maximal order of the
automorphism group for a genus 3 curve is 164 (given by the Hurwitz
bound), and $|\mathrm{SL}_2(\Z/7\Z)/\pm1|=164$, and the maximal
order of this group for genus 5 is 192 and
$|\mathrm{SL}_2(\Z/8\Z)/\pm1|=192$. The first lemma is proved by
Elkies in \cite{Elkies2}.

\begin{lema}\label{X7} The curve $X(7)$ is a genus 3 curve isomorphic over $\Q$ to the Klein
quartic which is a bielliptic curve and is not hyperelliptic.
\end{lema}

Recall that the Klein curve is the curve over $\Q$ defined by the
quartic equation $$x^3y+y^3z+z^3x=0.$$ Similarly, we take the model
$W$ defined over $\Q$ of the Wiman curve (which has the maximal
order of the automorphism group for a genus 5 curve) given as the
intersection of the following three quadrics in $\mathbb{P}^4$:
$$
x_0^2=x_3x_4 ,\ \ x_3^2=4x_1^2+x_2^2,\ \ x_4^2 = x_1x_2.
$$

\begin{lema}\label{X8} The curve $X(8)$ is a genus 5 curve isomorphic
over $\Q$ to the Wiman curve $W$, which is a bielliptic curve and is
not hyperelliptic.
\end{lema}

\begin{dem} One can easily see that $W$ is a curve with the same
group of automorphisms as $X(8)$ over $\C$. Since there is only
one such curve over $\C$, we get that they are isomorphic over
$\C$.

Consider the involution of $W$ over $\Q$ given by
$$\iota_1(x_0,x_1,x_2,x_3,x_4)=(x_0,x_1,x_2,-x_3,-x_4).$$
The quotient curve $W/\iota_1$ has equation given by
$$x_0^4=x_1x_2(4x_1^2+x_2^2)$$
which is isomorphic to the curve $X_0(64)$ over $\Q$ (e.g. a
computation via \cite{magma}).

Hence $X(8)$ and $W$ are curves over $\Q$, isomorphic over $\C$, and
both unramified degree covers of the same curve over $\Q$. Moreover,
one can see that there is only one involution of $W$ defined over
$\C$ whose quotient is $X_0(64)$: in fact, there are 4 involutions
without fixed points, three of them give quotients of genus three
that are hyperelliptic and one gives $X_0(64)$ (see also Subsection
3.2 in \cite{KMV}). We deduce that the cover $f:W\to X_0(64)$ must
be a twisted form (over $\Q$) of the cover $X(8)\to X_0(64)$.

The twisted forms of a fixed (degree 2) unramified covering are
well-known. In our case they can be described as the curves $W_d$
given by
$$
x_0^2=x_3x_4 ,\ \ dx_3^2=4x_1^2+x_2^2,\ \ x_4^2 = dx_1x_2,
$$
where $d$ is a squarefree integer, together with the natural map
$f$ to $X_0(64)$ given by
$f_d(x_0,x_1,x_2,x_3,x_4)=(x_0,x_1,x_2)$. Now, since the covering
$f$ is unramified, the subsets $f_d(W_d(\Q))$ do not intersect for
distinct covers and they give a partition of the set
$X_0(64)(\Q)$. This implies that only a finite number of covers do
have rational points. In our case, since $X_0(64)$ has 4
$\Q$-rational points corresponding to the cusps, a simple
computation shows that only $W=W_1$ and $W_2$ have $\Q$-rational
points; both have 4 $\Q$-rational points, each one covering 2
points of $X_0(64)(\Q)$.

Since $X(8)$ does have rational points (some cusps), we get that the
curve $X(8)$ is either isomorphic to $W$ or $W_2$ over $\Q$. But,
although $W$ and $W_2$ produce distinct coverings of $X_0(64)$ over
$\Q$, they are isomorphic as curves over $\Q$.
\end{dem}

\begin{obs2} A computation with MAGMA reveals that the curve $W$ is also
isomorphic over $\Q$ to the model of $X(8)$ over $\Q$ given by
Yang in \cite[Table page 507]{ya}. Moreover, it has a degree $2$
map to an elliptic curve isogenous to $X_0(32)$. In fact, $X(8)$
has three different bielliptic involutions (see \cite{KMV} or the
next section).
\end{obs2}

Finally, we recall that the curve $X_0(N^2)$ has genus $>1$ if and
only if $N\ge 8$. For them we have the following special case of a
result cited in the introduction.

\begin{prop}[Ogg-Bars, \cite{Og} \cite{ba}]\label{orba} The curve $X_0(N^2)$ with $g_{X_0(N^2)}\geq 2$ is never
hyperelliptic, and it is bielliptic exactly for $N=8$ and $9$.
\end{prop}

Now we can proceed to the proof of the main theorems in the
section.

\begin{dem1}
First of all, recall that, if $f:C\rightarrow C'$ is a
non-constant morphism between non-singular projective curves such
that the genus of $C'$ is $\geq 2$, and $C$ is hyperelliptic, then
$C'$ is hyperelliptic.

Since $X_0(N^2)$ is never hyperelliptic if the genus is larger than
1 by Proposition \ref{orba}, i.e. if $N>7$, we get that $X(N)$
cannot be hyperelliptic unless $N=7$. But this case was already
considered in Lemma \ref{X7}.
\end{dem1}

\begin{dem2}
Recall the following result by Harris and Silvermann in
\cite{HaSi}: Let $\phi:C\rightarrow C'$ be a non-constant morphism
between non-singular projective curves such that the genus of $C'$
is $\geq 2$. If $C$ is bielliptic then $C'$ is bielliptic or
hyperelliptic.

Now, the result follows for $N>9$ by using the map to $X_0(N^2)$
given in Lemma \ref{X0N2} and the result in Proposition \ref{orba}.
The cases $N=7$ and $8$ were already considered in Lemmas \ref{X7}
and \ref{X8}. It only remains to show that $X(9)$ is not bielliptic.

Recall that the genus of $X(9)$ is $10>6$. We will construct a map
$\rho$, which is a Galois cover (over $\C$) and it verifies the
conditions of the part (3) of Lemma \ref{lema3}. Consider the
(singular) model of $X(9)$ given by $y^6-x(x^3+1)y^3=x^5(x^3+1)^2$
\cite[p.507]{ya} (although this model is defined over $\Q$, we do
not know if it is isomorphic to $X(9)$ over $\Q$). Now let $E'$ be
the curve given by the equation $$z^2-x(x^3+1)z=x^5(x^3+1)^2.$$ We
get a map $\rho$ from $X(9)$ to $E'$ by taking $z=y^3$, which has
degree 3 (hence odd) and is Galois. The curve $E'$ is an elliptic
curve isomorphic to $E:t^2-t=x^3$ by writing
$t:=\frac{z}{x(x^3+1)}$. By applying part (3) of Lemma
\ref{lema3}, and since
\[2g(X(9))-2=2\cdot
10-2=18>\deg(\rho)(2g(E)+2)=3\cdot (2\cdot 1 +2)=12, \] we get
that $X(9)$ is not bielliptic.\end{dem2}

\begin{obs2} It is possible to describe theoretically the construction
in the last proof for the case $X(9)$.  First, consider the map
$f_0: X(9)\to X_0(81)$ given by Lemma \ref{X0N2}; it is a degree $3$
map to a curve of genus $4$. Then, consider the degree 3 map $\pi:
X_0(81)\to X_0(27)$, where the target is an elliptic curve. Finally,
let $E$ be the elliptic curve, given by the simple equation
$y^2-y=x^3$. The curve $E$ is 3-isogenous (over $\Q$) to the curve
$X_0(27)$. The map $\rho$ makes the following diagram commutative
$$ \xymatrix@R=0.7pc@C=0.8pc{
X(9) \ar@{->}[dd]^{\rho} \ar@{->}[rr]^{f_1} & & \ar@{->}[dd]^{\pi} X_0(81) \\
& & \\
E \ar@{->}[rr] & & X_0(27) }
$$
An analogous construction (but with degree 2 maps) can also be
done for the curve $X(8)$.
\end{obs2}
\end{section}

\begin{section}{On quadratic points for $X(N)$}

Let $C$ be a non-singular curve of genus greater than one, defined
over a number field $K$. Mordell's conjecture, proved by Faltings,
states that the set of $K$-rational points $C(K)$ of $C$ is always
finite. In order to generalize this, it is natural to consider the
set
$$\Gamma_d (C, K)=\bigcup_{[L : K]\leq d} C(L) $$ of points of
degree $d$ of $C$ over $K$. For quadratic points, that is, $d=2$,
Abramovich and Harris showed in \cite{abha} that $\Gamma_2(C,F)$ is
not finite for some finite extension $F$ of $K$ if and only if the
curve $C$ is either hyperelliptic or bielliptic. Hence, the
following result is a direct consequence of Theorems \ref{hypxn} and
\ref{biexn} in Section 4.

\begin{cor} The only modular curves $X(N)$ of genus $\geq 2$ such
that there exists a number field $L$ where the set
$\Gamma_2(X(N),L)$ is not finite are $X(7)$ and $X(8)$.
\end{cor}

Now, we can ask if, for $N=7$ or $N=8$, there are infinitely many
quadratic points over the cyclotomic field $\Q(\zeta_N)$ (which is
the smallest field where they can have non-cuspidal rational
points).

\begin{teo} For all $N\geq 7$, the number of quadratic points of $X(N)$ over
$F:=\Q(\zeta_N)$ is always finite.

\end{teo}
\begin{dem} By the corollary above, we
only need to study $N=7$ or $N=8$.

If $C(F)\neq\emptyset$, then by \cite{abha} we have:
$\#\Gamma_2(C,F)=\infty$ with $C$ a non-singular curve over $F$ if
and only if $C$ is hyperelliptic or has a degree two morphism
$\varphi:C\rightarrow E$ all defined over $F$,  with $E$ an
elliptic curve of $\operatorname{rank}_{\Z}E(F)\geq 1$.

It is known (see for example \cite{Pr}) that $\mathrm{Jac}(X(7))$
over $\Q(\zeta_7)$ is isomorphic to $E^3$, where $E$ is the elliptic
curve $y^2+3xy+y=x^3-2x-3$, which is isomorphic to $X_0(7^2)$.
Therefore, since $X(7)$ is non-hyperelliptic, we have an infinite
number of quadratic points over $\Q(\zeta_7)$ only if there is a
degree two map $X(7)\rightarrow E'$ all defined over $\Q(\zeta_7)$,
where $E'$ is an elliptic curve of positive rank over $\Q(\zeta_7)$.
But then, necessarily, $E'$ is isogenous to $E$ and, in particular,
${\rank_{\Z}}E'(\Q(\zeta_7))={\rank_{\Z}}(E({\Q(\zeta_7)})$. But
this last rank is zero, as a (2-Selmer) computation with
\verb+MAGMA+ \cite{magma} or \verb+SAGE+ \cite{sage} reveals.

For $N=8$, consider the equations over $\Q$ given above. Some
computations with \verb+MAGMA+ shows that the group of
automorphisms over $\Q$ is abelian and isomorphic to $(\Z/2\Z)^3$,
and the quotient with respect two of the elements gives the
elliptic curve $E$ with equation $y^2=x^3-x$ of conductor 32, and
by a third element gives the elliptic curve $E'$ with equation
$y^2=x^3+x$, of conductor 64. By \cite{KMV} there are exactly 3
bielliptic involutions for $X(8)$, so these are all of them. The
elliptic curves $E$ and $E'$ become isomorphic over $\Q(\zeta_8)$.
Hence, they have the same rank. Finally, a (2-Selmer) computation
with \verb+MAGMA+ or \verb+SAGE+ reveals that
${\rank_{\Z}}(E(\Q(\zeta_8)))=0$, proving the result.
\end{dem}

\begin{obs2} A computation with the help of the \verb+MAGMA+ algebra system
\cite{magma} shows that
$$\#\Gamma_2(X(8),\Q(\zeta_8))= 24,$$ corresponding to the cusps. This
result is done through computing all the quadratic points of
$X_0(8^2)$ over $\Q(\zeta_8)$. The curve $X_0(8^2)$  is a genus
$3$ curve, with Jacobian isogenous to the cube of the elliptic
curve $X_0(32)$ over $\Q(\zeta_8)$, which has only a finite number
of points over $\Q(\zeta_8)$. Then we compute the inverse image
with respect to the the degree $2$ map from $X(8)\rightarrow
X_0(32)$. We obtain that although there are points in
$\Gamma_2(X_0(32),\Q(\zeta_8))$ which do not come from cusps
(there are more than one hundred points), none of them lift to a
quadratic point of $X(8)$.
%They analogous result is also true for $X(7)$.
%See the second author homepage for the MAGMA files.
\end{obs2}

\end{section}

\begin{section}*{Acknowledgements}
 The second author wishes to thank G. Cornelissen for pointing him
to the problem of automorphisms of $X(N)$. The first and third
author would like to thank J. Gonz\'alez for warning us about some
results by Carayol and for his useful remarks on the writing of
this note, also they thank E. Gonz\'alez-Jim\'enez for his
comments. We are pleased to thank A. Schweizer for his suggestions
and comments. We thank our colleague H.D. Stainsby by
proof-reading our English. Finally, we thank an anonymous referee
for all the remarks and suggestions in order to improve the paper.
 \end{section}


\begin{thebibliography}{70}
\bibitem{abha}D. Abramovich and J. Harris, Abelian varieties and
curves in $W_d(C)$, Compositio Math. 78 (1991), 227-238.

\bibitem{AkSi} M. Akbas and D. Singerman, The normalizer of $\Gamma_0(N)$
in $\operatorname{PSL}_2(\mathbb{R})$, Glassgow Math. J. 32 (1990),
317-327.

\bibitem{ba} F.Bars, Bielliptic modular curves, J. Number Theory 76 (1999),
154-165.

%\bibitem{ba2} F. Bars, The group structure of the normalizer of $\Gamma_0(N)$ after Atkin-Lehner. Communications in Algebra 36, No. 6, 2160-2170
%(2008).

%\bibitem{cre} Cremona tables. See http://homepages.warwick.ac.uk/~masgaj/ftp/data/count.00000-09999.gz

\bibitem{DR} P. Deligne, M. Rapoport, Les sch{\'e}mas de modules des courbes
elliptiques in Modular Functions of One Variable II, Springer
Lecture Notes in Mathematics 349 (1973), 143-316.

\bibitem{Gur} P.Bendig, A.Camina and R.Guralnik, Automorphisms of
the modular curve. Progress in Galois Theory, Dev.Math. vol.12,
25-37, Springer (2005).


\bibitem{magma} Wieb Bosma, John Cannon, and Catherine Playoust, The Magma algebra system. I. The user language,
J. Symbolic Comput. 24 (1997), 235-265.


\bibitem{CKK} G.Cornelissen, F. Kato and A. Kontogeorgis, Discontinuous groups in positive characteristic and
              automorphisms of {M}umford curves. Math. Ann.320,
              No.1, (2001), 55-85.


\bibitem{Elkies} N.D. Elkies, The automorphism group of the modular curve $X_0(63)$, Compositio Math. 74 (1990),
203-208.

\bibitem{Elkies2} N.D. Elkies, The Klein quartic in number theory, pages 51-102 in The Eightfold
Way: The Beauty of Klein's Quartic Curve, S.Levy, ed.; Cambridge
Univ. Press, 1999.

\bibitem{F-K} H.M. Farkas and I.Kra, Riemann surfaces. Graduate
Texts in Mathematics vol. 71, second edition. Springer (1992).

\bibitem{HK} E. Halberstadt and A. Kraus, Sur la courbe modulaire $X_E(7)$,
Experiment. Math. 12 (2003), no. 1, 27-40.


\bibitem{HaSi} J.Harris and J.H.Silverman, Bielliptic curves and
symmetric products. Proc. Am. Math. Soc. 112 (1991), 347-356.

\bibitem{Har} M. Harrison, A new automorphism of $X_0(108)$,
http://arxiv.org/abs/1108.5595, (2011).


\bibitem{IsMo} N.Ishii and F. Momose, Hyperelliptic modular curves,
Tsukuba J. Math. 15 (1991), 413-423.

\bibitem{jeki} D. Jeon and C.H.Kim, Bielliptic modular curves
$X_1(N)$, Acta Arith. 112.1 (2004), 75-86.

\bibitem{jeki2} D. Jeon and C.H.Kim, Bielliptic modular curves $X_1(M,N)$, Manuscripta
Math. 118 (2005), 455-466.
%
\bibitem{jeki3} D. Jeon, C.H. Kim and A.Schweizer, Bielliptic intermediate
modular curve, work in progress.

\bibitem{KMV} T.Kato, M.Magaard and H.V\"olklein, Bi-elliptic Weierstrass points on curves of genus
5. Indag.Math. (N.S.) 22, (2011) no.1-2, 116-130.

\bibitem{KM} N. M. Katz B. and Mazur, Arithmetic Moduli of Elliptic
Curves. Annals of Mathematics Studies 108, Princeton University
Press (1985).

\bibitem{KeMo} M.A. Kenku and F. Momose, Automorphsism groups of
modular curves $X_0(N)$. Compositio Math. 65 (1988), 51-80.

\bibitem{kiko} C.H. Kim and J.K. Koo, The normalizer of
$\Gamma_1(N)$ in $\operatorname{PSL}_2(\R)$, Comm. Algebra 28
(2000), 5303-5310.

%\bibitem{ko} A. Kontogeorgis, On automorphisms of certain algebraic curves and varieties,
%Ph.D. dissertation, University of Crete 1999 (In Greek)

\bibitem{Mes} J-F. Mestre, Corps euclidiens, unit\'es
exceptionnelles et courbes \'elliptiques. J.Number Theory 13
(1981), 123-137 .


\bibitem{Mo} F. Momose, Automorphism groups of the modular curves
$X_{1}(N)$. Preprint.

%\bibitem{Ne} M. Newman, Struture theorem for modular subgroups. Duke
%Math. J. 22, 25-32 (1955).

\bibitem{Og} A.P. Ogg, Hyperelliptic modular curves. Bull.Soc.Math.
France 102 (1974), 449-462.

\bibitem{Pr} D.T.Prapavessi, On the Jacobian of the Klein curve.
Proc. Amer. Math. Soc., 122 (4) (1994), 971-978.

\bibitem{Ritz} C.Ritzenthaler, Automorphismes des courbes modulaires
$X(n)$ en caract\'eristique $p$. Mauscripta Math. vol. 109 (2002),
49-62.

\bibitem{se} J.-P. Serre, The automorphism group of X(p), London Math. Soc. Lect. Notes Series 254
(1998), appendix to B. Mazur, Open problems regarding rational
points on curves and varieties, in \textit{Galois representations
in arithmetic algebraic geometry}, A. J. Scholl and R. L. Taylor,
eds.

\bibitem{Sch} A. Schweizer, Bielliptic Drinfeld modular curves.
Asian J. Math. 5 (2001), 705-720.

\bibitem{Shi} G.Shimura, Introduction to the arithmetic theory of
automophic functions. Publ. of The Mathematical Society of Japan,
11, Princeton Univ.Press (1971).

\bibitem{sage} {\sc W. Stein et al.}, Sage: {O}pen {S}ource {M}athematical {S}oftware ({V}ersion 4.3), The  Sage~Group, 2009,
 {\tt http://www.sagemath.org}.


\bibitem{ya} Y.Yang, Defining equation of modular curves, Advances
in Mathematics 204 (2006), 481-508.
\end{thebibliography}
\end{document}